\documentstyle[12pt]{article}
\def\CC{{\rm\kern.24em\vrule width.02em height1.4ex depth-.05ex\kern-.26em C}}
\def\BB{{\rm\kern.24em\vrule width.02em height1.4ex depth-.05ex\kern-.26em B}}
\def\HH{{\rm\kern.24em\vrule width.02em height1.4ex depth-.05ex\kern-.26em H}}
\def\PP{{\rm\kern.24em\vrule width.02em height1.4ex depth-.05ex\kern-.26em P}}
\def\RR{{\rm\kern.24em\vrule width.02em height1.4ex depth-.05ex\kern-.26em R}}

\newtheorem{prob}{Problem }[subsection]
\def\bp{\begin{prob}}
\def\ep{\end{prob}}
\newtheorem{oprob}{Prize Problem }
\def\bop{\begin{oprob}}
\def\eop{\end{oprob}}

\begin{document}

\title{ {\bf THE SEVERAL COMPLEX}\\ {\bf VARIABLES PROBLEM LIST}\\
{\bf IT'S BACK!!!!!!!!}}
\author{ }
\date{September 1995}
\maketitle

During the past two years, we got many requests from people
interested in the list, but we did NOT get many new problems.
Since we believe that this list is worth pursuing, we are
making a new attempt.  The SCV problem list is an evolving, 
growing part of our
subject. Please contribute your ideas for helping to keep
it lively. In this come--back edition we request from you:

\begin{itemize}
\item[1.-] {\bf NEW PROBLEMS!!}
\item[2.-] Tell us if you know that any of the problems
have been recently solved. We will announce it in the
next edition.
\end{itemize}

To simplify the reading, we are stating the new problems for this
edition at the beginning of the list, right after the introduction.

\bigskip

\section{ Introduction}

The purpose of this bulletin board is to collect problems in higher
dimensional complex analysis. We are interested both in basic research
questions as well as interactive questions with other fields and sciences.  
We encourage everybody to submit problems to the list. This 
includes not only those coming up in your own work, but also
others- maybe well known and classical- that you see missing, but that
you think workers in the field should be aware of. 
Not only are we searching for basic research type questions in several
complex variables, we also solicit questions exploring relations to 
other mathematical fields, one complex variable, partial differential
equations, differential
geometry, dynamics, etc. and to other sciences such as physics,
engineering, biology etc. 
While some questions fall rather naturally into one
of the subject areas in this problem list and may
lead to a publishable paper, other questions may be non-
specific or of a transient or technical quality. Thus we have
a section called ``Scratchpad'' for conversational questions,
vaguely formulated questions, or questions to which you may hope
to get a quick answer. In the ``open prize problems'' section,
 you are welcome to offer a
nice little prize to whoever does it. The 
miscellaneous section is for announcements of conferences, jokes, remarks
on the general state of the field, etc.  

This file is in Latex. ( In order to get a table of contents you 
should run the file twice.) If a problem appears in the problem list
the first time
or has been changed, it will be marked by the word {\bf{New}}.  

\bigskip

{\em To submit problems and other suggestions, email to:}  

\begin{itemize}
\item  Gregery Buzzard at gbuzzard@iu-math.math.indiana.edu, 
\item  John Erik Fornaess at fornaess@math.lsa.umich.edu , 
\item  Estela A. Gavosto at gavosto@math.lsa.umich.edu or to  
\item  Steven G. Krantz at sk@math.wustl.edu
\end{itemize}

\noindent {\bf ANNOUNCEMENTS: }

\begin{itemize}
\item   Academic Year 95/96, Special Year in Several Complex Variables,
 MSRI, Berkeley. 
\item January, 10-13, 1996 Orlando, FL (1996 Joint Meetings)
\begin{itemize}
\item AMS Special Session on Analytic Methods in SCV.
\item AMS Special Session on Multidimensional Complex Dynamics.
\end{itemize}
\item March 22-23, 1996 Iowa City, IA . AMS Special Session on
Geometric and Analytic Methods in SCV
\end{itemize}

\newpage

\section{NEW PROBLEMS IN THIS EDITION }

\bp {\bf (New) (See Section 7.2.) Let $V$ be a germ at $0$ of an irreducible complex variety
 in $\CC^n.$ For any small enough $\epsilon > 0$, let
$E= E_\epsilon:= V \cap \RR^n \cap \{|z|< \epsilon \}.$
Define the extremal function $U_E$ on $\{ |z|< \epsilon \}:
U_E(z):= \sup \{u(z); u \leq 0\; \mbox{on}\; E,\;
0 \leq u \leq 1 \; \mbox{on} \; \{|z|< \epsilon \}$
The problem is to classify those $V$ for which we have a constant
$A>0$ for which
$(\alpha) \; U_E(z) \leq A | \Im z |,\; |z|< \epsilon/2.$}
\ep

\bp {\bf (New) (See Section 8.2.) What is the lowest degree of a polynomial mapping
$P: \RR^2 \rightarrow \RR^2$ for which the (real) Jacobian conjecture
fails, i.e. the Jacobian of $P$ near vanishes while $p$ fails
to be invertible.}
\ep

\bp {\bf (New) (See Section 11.3) (Liouville type assertion on CR 3-manifolds) On a CR 3-manifold
one can always define a Carnot-Caratheodory metric. A homeomorphism between
two strongly pseudoconvex CR 3-manifolds is called conformal if it maps
infinitesimal spheres with respect to a Carnot-Caratheodory metric to
infinitesimal spheres. Is a conformal homeomorphism CR?  
-Puqi Tang, tang@math.purdue.edu}
\ep
\bp  {\bf (New) (See Section 11.3) A diffeomorphism between two strongly pseudoconvex CR
$(2n+1)$-manifolds is called quasiconformal if its differential preserves
the underlying contact structures and distorts the CR structures boundedly.
When Hermitian metrics on the contact bundles are fixed, this distortion
can be measured by checking how spheres in the contact space is mapped to
ellipsoids. However, this measurement depends on the choices of the Hermitian
metrics if $n>1$. Fixing a quasiconformal diffeomorphism, can we choose 
Hermitian metrics so that the distortion is minimal? 
- Puqi Tang, tang@math.purdue.edu}
\ep

\newpage

\tableofcontents

\newpage

\section{ Envelopes of holomorphy of domains in $\CC^n$}
\subsection{ C-R functions}
\bp Let M be a smooth real hypersurface in $\CC^n$ with defining
function r. Let $p_0$ be a point in M, and let U be a neighborhood in $\CC^n$
of $p_0$. When do C-R functions on $M\cap U$ extend to $\lbrace r > 0
\rbrace$?\\
-B. Stensones, berit@math.lsa.umich.edu
\ep
\subsection{ The Future Tube}
\subsection{ Schlichtness of envelopes of holomorphy}
\bp Let D be a bounded domain in $\CC^n$ with smooth boundary, with
envelope of holomorphy $D^{\sim}$. 
A)Is  $D^{\sim}$ finitely sheeted over $\CC^n$?
B) Does $D^{\sim}$ have finite volume?\\
 -B. Stensones, berit@math.lsa.umich.edu.
\ep
\bp
If D is a strictly pseudoconvex domain in $\CC^2$ diffeomorphic to $\RR^4$ and A is an embedded
analytic disc outside D with $bA \subset bD$, then there is an embedded smooth disc
M in D, with bM=bA, and such that M is totally real except at
one elliptic point. With such a configuration
(M,A) in $\CC^2$, does it follow that the envelope of holomorphy of
M contains A ?\\
-F. Forstneric, forstner@math.wisc.edu
\ep

\newpage

\section{ The Levi Problem}
\subsection{  Locally Stein open subsets of Stein Spaces}
\bp Let $U$ be an open subset of a Stein space. If $U$ is locally Stein, 
is $U$ Stein?\\
\ep  

\bp Let M be a compact complex manifold. Find some appropriate
conditions on M such that if U is an open subset of M which is 
locally Stein, then U is Stein.\\
-N. Sibony, sibony@anh.matups.fr
\ep
\subsection{ The Runge Problem}
\bp Let $X$ be a complex subvariety of the unit ball $\BB$ with an isolated
singularity at the origin. Suppose $\Omega_t$ is a continuously increasing
family of Stein open subsets of $X$. Also, suppose that $0 \in \partial \Omega_{1/2} ,
0\in \Omega_t, t>1/2.$ Is $\Omega_s$ Runge in $\Omega_t$ $\forall s<t$?
\ep
\bp Let $K^{compact} \subset X \backslash \lbrace 0 \rbrace$, $X$ as above, 
$0 \in \hat{K}$, the polynomially convex hull of $K$. Suppose that $L \subset X$,
$L$ contains a neighborhood of $K$ in $X$. Does $\hat{L}$ contain a neighborhood
of $0$ in $X$? (This is equivalent to the previous problem.)\\

\ep
\subsection{ The Union Problem}
\bp Let $\lbrace U_n \rbrace$ be an increasing sequence of Stein open subsets of
a Stein Space $X$. Is $\bigcup_{n=1}^{\infty} U_n$ Stein? (If X is a complex
manifold, not necessarily Stein, the answer is in general no.)\\

\ep
\bp Suppose $M$ is a complex manifold of dimension 2, and suppose $\forall K^{compact}
\subset M$ there exists an open subset $U$, $K \subset U \subset M$ such that 
$U$ is biholomorphic to the unit ball. Does this imply that M is equivalent to
the unit ball $\BB$, $\Delta \times \CC$ or $\CC^2$? (This is false in $\CC^3$).\\

\ep
\bp Is ``long'' $\CC^2$ biholomorphic to $\CC^2$? ( A complex manifold is a long
$\CC^2$ if it is the union of proper subsets that are biholomorphic to $\CC^2$.)\\
\ep

\newpage

\section{Holomorphic function theory of domains in $\CC^n$}
\subsection{ Approximation Problems}
\bp Let $\Omega \subset \subset \CC^n$ be a smooth pseudoconvex 
domain. Let $A^k(\bar\Omega)= H(\Omega)\cap C^k(\bar\Omega)$
with $C^k$ topology, $1  \leq k \leq \infty$. Is $A^{\infty}(\bar\Omega)$
dense in $A^k(\bar\Omega)$?\\
-N. Sibony, sibony@anh.matups.fr
\ep
\bp Let $\Omega \subset \subset \CC^n$ be a smooth pseudoconvex 
domain. Assume that $\bar\Omega$ has a Stein neighborhood basis. Is 
every function in $A^{\infty}(\bar\Omega)$ uniformly approximable
by holomorphic functions in a neighborhood of $\bar\Omega$?\\
-N. Sibony, sibony@anh.matups.fr
\ep

\subsection{ The corona problem for the ball or the polydisc}     

\bp
 Let $f_1, ..., f_k$ be bounded holomorphic functions on $\Omega.$
Suppose $\sum |f_i|> c> 0.$ The Corona problem is whether there exist
bounded holomorphic funcions $g_1, ..., g_k$ on $\Omega$ such
that $\sum f_i g_i \equiv 1.$ Can the Corona problem be solved if $\Omega$
is the unit ball or the unit polydisc?
\ep

\subsection{Hyper-holomorphic functions}
\bp  \ \ \
Let $\Phi =\{ \Phi^0, \Phi^1, \Phi^2, \Phi^3\}$ be a subset of the body 
$(=skew-field)  \ \HH$ of the real quaternions with the conditions
$$
\Phi^p \bar\Phi^q + \Phi^q \bar\Phi^p = 2 \delta_{p,q},
$$
$\Omega$ is a domain in $\RR^4=\CC^2$. On the $\HH$-bimodule $C^1(\Omega ;
\HH)$ the operator
      $D_\Phi = \sum^3_{k=0} \Phi^k \frac \partial {\partial x_k} $
defines the set of the $(left-\Phi)$-hyper-holomorphic functions 
$M_\Phi (\Omega ; \HH) = Ker D_\Phi$.
Let $\hat \Phi =\{ 1, i, j, -k\}$ with i, j, k being the ordinary imaginary
units in $\HH$. It is known that
$$
Hol (\Omega ; \CC^2) \subset M_{\hat \Phi} (\Omega ; \HH),
$$
where $Hol (\Omega ; \CC^2)$ is the subset of all holomorphic (in the usual
sense) mappings.\par
1) Describe the set of all hyper-holomorphic but non-holomorphic functions in 
$\Omega$ (it is not empty).\par
2) Describe "locations" of the sets $Hol$ and $M_\Phi Hol$ inside $M_\Phi$.\par
3) How do they depend on $\Omega$?

-H. Shapiro and N. Vasilievsky,  matemat@cinvesmx.bitnet
\ep

\subsection{ Peak Points and Peak Sets}
\bp Is every boundary point of a bounded pseudoconvex domain $D$ in $\CC^n$ of
finite type a peak point for $A(D)$?
\ep
\bp Describe peak sets on weakly pseudoconvex domains of finite type.
\ep
\subsection{ Representing measures and polynomial hulls}
\bp  Characterize the representing measures for the origin in
the unit ball in $\CC^n.$
\ep 

\bp
  Let $X \subset \CC^n$ be compact
and totally disconnected.  Assume that $X$ has finite one-dimensional Hausdorff
 measure.  Then is $X$ polynomially convex?\\
-H. Alexander, U22330@UICVM.bitnet
\ep
\subsection{Zero sets of holomorphic functions}
\bp Let $\Omega \subset \CC^n$ be a smoothly bounded, pseudoconvex domain.
Let $V \subseteq \Omega$ be a divisor which is bounded away from the
Silov boundary in $\partial \Omega.$  Prove that there is a bounded
holomorphic function on $\Omega$ that vanishes on $V.$  (This is known
for the polydisc---a result of Rudin.)\\ 
-S. Krantz, sk @math.wustl.edu
\ep
\bp Give a geometric characterization of zero sets of $H^p$ functions on 
strongly pseudoconvex domains.
\ep

\newpage

\section{ Existence and regularity properties of the Cauchy Riemann
operator}
\subsection{ The $\bar\partial$ Neumann problem}
\bp Solve the $\bar\partial\; Neumann$ problem in $H^s$.
\ep
\bp Solve $\bar\partial \; Neumann$ in $L^2 (d\mu)$, $d\mu$ a measure, even on the
unit disc.
\ep

\bp 
Let $\Omega$ be a smooth bounded planar domain and let $T:L^2(\Omega)\mapsto
L^2(\Omega)$ be the Kohn solution operator for $\partial/\partial\bar z$.  For
$h\in L^\infty(\Omega)$ we define the compact non-self-adjoint operator
$S_h:L^2(\Omega)\mapsto L^2(\Omega), f\mapsto \bar hTf$.  Under what circumstance
will $S_h$ have eigenvalues?   Consider in particular the case where $h$ is
holomorphic.  If $h$ has a single-valued holomorphic primitive $H$ then it is
not hard to see that $S_h$ has no eigenvalues.  (Just note that $e^{-\bar
H/\lambda}Tf$ is holomorphic and orthogonal to holomorphic functions.)
Does the converse hold?  That is:  do periods of $H$ guarantee the existence
of eigenvalues?  (For annuli, the answer is``yes.")\\
-D. Barrett, barrett@math.lsa.umich.edu
\ep

\bp
 If D is a real analytically bounded pseudoconvex domain in $\CC^n$,
 is the $\bar\partial$ Neumann problem globally real analytic hypoelliptic 
up to the boundary?\\
- D. S. Tartakoff, U22393@UICVM.bitnet
\ep
\bp 
If M is a compact real analytic pseudoconvex CR manifold
          of real dimension at least 5 is $\Box_b$ globally real
 analytic hypoelliptic (say on (0,1) forms) on M?\\
- D. S. Tartakoff, U22393@UICVM.bitnet
\ep
\subsection{ H\"{o}lder estimates for $\bar \partial$ on pseudoconvex
domains  of finite type in $\CC^n$}

\bp Can one solve $\bar\partial$ with  H\"{o}lder estimates on any real 
analytic bounded pseudoconvex domain in $\CC^n$?\\
\ep

\subsection{ $L^p$ Estimates for $\bar\partial$}
\bp Solve $\bar\partial$ in $L^p$ in bounded pseudoconvex domains, without necessarily
 smooth boundary, $1 < p < 2$. 
\ep
\bp Solve $\bar\partial$ with $L^p$ estimates, $1< p< \infty$ for smooth convex
 domains in $\CC^n$.
\ep
\subsection{Uniform estimates for $\bar \partial $ on pseudoconvex domains}
\bp Can one solve $\bar \partial u = f$ in sup norm on smooth bounded convex
 domains in $\CC^2$?\\
\ep
\bp Which geometric properties of the domain do uniform estimates depend on?
 $\lbrack e.g. L^2$ estimates depend on the diameter and the dimension
 $\rbrack$.
\ep
\bp If $\Omega \subset \subset \CC^n$ is pseudoconvex with smooth boundary, does
 $\bar \partial: \Lambda^{\alpha}(\bar \Omega) \mapsto \Lambda_{(0, 1)}^{\beta}
 (\bar \Omega)$ have closed range if $\alpha < \beta$?
\ep

\bp
Let $\Omega$ in $\CC^n$ be a bounded biholomorphic image of the unit ball
in $\CC^n$. Does the $\bar\partial$ satisfy uniform estimates on $\Omega$?
The problem is of interest mainly when $\Omega$ has very rough boundary.\\
-B. Berndtsson, bob@math.chalmers.se
\ep

\newpage

\section{Geometric and Topological Properties of Pseudoconvex Domains}
\subsection{ Topological Properties of Pseudoconvex Domains}
\bp 
 Does there exist a pseudoconvex domain in $\CC^2 $ with the same homotopy
 type as $S^2$?\\
-F. Forstneric, forstner@math.wisc.edu
\ep

\subsection{ The Bergman Metric B}
\bp Is B complete on bounded pseudoconvex domains with less than $C^1$ 
boundary?\\
-J. McNeal, mcneal@math.princeton.edu
\ep
\bp For which domains is the Bergman kernal function $K(z, w)$ nonvanishing 
everywhere? This is related to Lu Qi-Keng conjecture.
\ep

\subsection{ The Caratheodory Metric c}
\bp Is c complete on a bounded smooth pseudoconvex domain?
\ep
\bp Is there a smoothly bounded pseudoconvex domain in $\CC^n$ on which the
 Kobayashi and Caratheodory metrics are not comparable?
 
-S. Krantz, SK@math.wustl.edu
\ep

\subsection{ Finite Type}
\bp Is the regular order of contact of a smooth pseudoconvex domain the same as
 the commutator type, one vector field at a time.\\
-J. McNeal, mcneal@math.princeton.edu
\ep
\bp Are finite type pseudoconvex domains in $\CC^2$ locally convexifiable 
with biholomorphic
 maps with continuous boundary values?
\ep

\subsection{ The Kobayashi Metric k}
\bp Is k complete on a bounded smooth pseudoconvex domain?
\ep
\bp What's the asymptotic behavior of k on finite type domains in $\CC^n$?
\ep
\bp
Let $D = \lbrace z \in \CC^n ; h(z) < 1 \rbrace$ be a bounded pseudoconvex 
complete$-$circular domain with
Minkowski function h. Characterize the Caratheodory, the Bergman or the
Kobayashi completeness of D via properties of h.
[The following is known: a) If h is continuous, then D is Bergman
complete.  b) If $n\geq 3$, continuity of h does not imply the Kobayashi
completeness of D.]\\
-P. Pflug, pflugvec@dosuni1.bitnet
\ep
\bp 
Given a compact complex hypersurface X in $\PP^n$ with degree at least $n+2$ and
only normal crossing singularities. Show that the infinitesimal Kobayashi
 pseudometric
on $\PP^n \setminus X$ is degenerate at most on a complex hypersurface.\\
-N. Sibony, sibony@anh.matups.fr
\ep
\bp
Find (or prove it does not exist) a complete hyperbolic manifold not 
biholomorphic to a 
convex domain where the Kobayashi and Carath\'eodory distances agree 
everywhere.\\
-M. Abate, abate@vaxsns.infn.it
\ep
\bp 
Can you recover certain
derivatives of the Levi form at a weakly pseudoconvex point in terms
of the weighted boundary behavior of the Kobayashi metric?\\
-S. Krantz,  sk@math.wustl.edu
\ep

\subsection{ Regular Domains}
\bp Let $\Omega \subset \subset \CC^n$ be pseudoconvex with smooth boundary.
 For $K \subset \bar\Omega$ compact, let $\hat K:= \lbrace z \in \bar \Omega ;
 \phi(z) \leq sup_K \phi, \phi$ continuous on $\bar \Omega, \phi$ 
plurisubharmonic on $\Omega \rbrace$. Assume $\hat K \cap \partial \Omega = K
 \cap \partial \Omega$.
 Then $\bar \Omega$ has a Stein neighborhood basis.\\
-N. Sibony, sibony@anh.matups.fr
\ep

\newpage

\section{ Plurisubharmonic Functions}
\subsection{The Complex Monge Ampere Operator}
\bp Find conditions insuring that if a sequence of plurisubharmonic functions 
$u_m \rightarrow u$ in $\CC^n$ then $(dd^c u_m)^n \rightarrow (dd^c u)^n$.
\ep
\subsection{Extremal Functions}
\bp {\bf (New) Let $V$ be a germ at $0$ of an irreducible complex variety
 in $\CC^n.$ For any small enough $\epsilon > 0$, let
$E= E_\epsilon:= V \cap \RR^n \cap \{|z|< \epsilon \}.$
Define the extremal function $U_E$ on $\{ |z|< \epsilon \}:
U_E(z):= \sup \{u(z); u \leq 0\; \mbox{on}\; E,\;
0 \leq u \leq 1 \; \mbox{on} \; \{|z|< \epsilon \}$
The problem is to classify those $V$ for which we have a constant
$A>0$ for which
$(\alpha) \; U_E(z) \leq A | \Im z |,\; |z|< \epsilon/2.$}
\ep

\newpage

\section{ Holomorphic maps}
\subsection { Automorphisms}
\bp Which smooth domains have non compact automorphism groups?\\
\ep
\subsection {The Jacobian Conjecture}
\bp 
If P: $\CC^2 \rightarrow \CC^2$ is a holomorphic polynomial map with nonvanishing
Jacobian, then P is biholomorphic.
\ep
\bp {\bf (New) What is the lowest degree of a polynomial mapping
$P: \RR^2 \rightarrow \RR^2$ for which the (real) Jacobian conjecture
fails, i.e. the Jacobian of $P$ never vanishes while $p$ fails
to be invertible.}
\ep

\subsection{Proper correspondences}
\bp
Classify Proper Holomorphic Correspondences from the unit disc to itself.
\ep
\bp 
 Is every proper holomorphic self-mapping
of a smooth bounded domain in $\CC^n (n>1)$ biholomorphic?
\ep

\subsection{ Real analytic domains (not necessarily pseudoconvex) in $\CC^n$, $f: U \rightarrow V$ biholomorphic}
\bp Does f have a continuous extension to $\bar U$? (This is true if U, V are
 pseudoconvex.)
\ep

\subsection {Stein manifolds}

\bp  If M is an n-dimensional parallellizable Stein manifold, can M
be holomorphically immersed in $\CC^n$?
[ True for n=1 (Gunning and Narasimhan); also such an M can be immersed in
 $\CC^{n+1}$ (Eliashberg-Gromov)].
\\
-R. Narasimhan, Dept. of Math., Univ. of Chicago, Chicago, IL 60637 
\ep
\bp Can $SL(n, \CC)$ be holomorphically immersed in $\CC^{n^2-1}$.
[Cannot be done algebraically, true if n=2]\\
-R. Narasimhan, Dept. of Math., Univ. of Chicago, Chicago, IL 60637 
\ep
\bp Can the algebraic hypersurface
$$ w^3 + w g(w_1, ......, w_n) + f(w_1, ......, w_n) = 1, $$
where g is homogeneous of degree 2 and f homogeneous of degree 3, be 
holomorphically immersed in $\CC^n$
[ Not true for algebraic immersion]\\
-R. Narasimhan, Dept. of Math., Univ. of Chicago, Chicago, IL 60637 
\ep
\bp  Suppose that the polynomial equation 
$$F(w_0, w_1, ......, w_n)=1$$
 defines a smooth
affine algebraic hypersurface (these are always parallellizable 
(Murthy, Swan, n=2;  Suslin for general n). Can they be holomorphically 
immersed in
$\CC^n$?\\
-R. Narasimhan, Dept. of Math., Univ. of Chicago, Chicago, IL 60637 
\ep

\newpage

\section{ Dynamical properties of holomorphic maps}
\subsection{ Dynamics of entire holomorphic maps}
\bp Do there exist Fatou--Bieberbach domains which are not Runge?
A Fatou--Bieberbach domain is a proper subdomain of $\CC^2$ which is 
biholomorphic to $\CC^2$.\\
-F. Forstneric, forstner@math.wisc.edu
\ep



\subsection{ The invariant set $K^+$ of a polynomial complex Henon map}
\bp Do there exist wandering domains for complex Henon maps?
\ep
\bp The Fatou Set of a Polynomial Automorphism of $\CC^2$.
Let $f$ be a polynomial automorphism of $\CC^2$, and suppose that $f$ is
not conjugate to an elementary map.  Let $\Omega$ be a connected component of
the set of points where the forward iterates $\lbrace f^n:n=1,2,3,...\rbrace $
are locally bounded, which is the Fatou set of $f$.  Let us suppose that 
$\Omega$
is periodic, i.e. $f^n\Omega=\Omega$ for some $n\ne0$.
Let us first consider the case in which $\Omega$ is recurrent, i.e. $\Omega$
contains a point whose orbit does not converge to $\partial\Omega$.  In
this case $\Omega$ is either (a) a basin of attraction of a sink orbit, (b) the
basin of attractin of a Siegel disk, or (c) the basin of a Herman ring.  It is
not hard to construct examples of (a) and (b).\\
\noindent(i)  Can case (c) occur?
Nothing at all seems to be known about periodic domains which are not 
recurrent.\\
\noindent(ii)  If $\Omega$ is not recurrent, does there exist
$P\in\partial\Omega$ such that $\lim_{k\to\infty}f^{nk}(Q)=P$ for all
$Q\in\Omega$?\\
\noindent(iii)  If $\Omega$ is not recurrent, does there exist
$P\in\partial\Omega$ such that $f^nP=P$, and one of the eigenvalues of $DF(P)$
is $e^{2\pi ip/q}$?\\
-E. Bedford, bedford@iubacs.bitnet or bedford@ucs.indiana.edu.
\ep

\bp
Let $ F_c : (z,w) \rightarrow (z^2 +c - w, z)$, $c \in \CC$ a symplectic
automorphism of $\CC^2$, i.e. $F_c$ is biholomorphic and $ F^*_c(dz \wedge dw) = dz \wedge dw$. Can there exist a $c \in \CC$ and a periodic orbit $\{z_i\}^k_{i=0}$ for $F_c$ such that $z_0$ belongs to a Siegel domain ?\\
-M. Herman
\ep

\subsection{ Dynamics on $\PP^k$}

\bp Suppose $f$ is a holomorphic map on $\PP^2$ of degree $d \geq 2$. 
Must $f$ have a repelling fixed point.\\
\ep

\bp
 Does every holomorphic map $F:\PP^k \rightarrow \PP^k$ of degree at least two have a repelling periodic point? Do meromorphic maps have periodic orbits?
-sibony@anh.matups.fr
\ep

\bp
 Is the support of $\mu$ equal to the closure of the repelling periodic orbits?
\ep

\bp
 Does there exist a holomorphic map $F:\PP^k \rightarrow \PP^k$ with a wandering Fatou component $U$, i.e. $F^n(U) \bigcap F^m(U) = \emptyset $ for all $n \neq m$?
\ep

\bp
Classify the dynamics around a  fixed point of a holomorphic map on $\PP^k$
or even just defined in a neighborhood of p.
-sibony@anh.matups.fr
\ep

\bp
Let $f:\PP^2 \rightarrow \PP^2$ be a holomorphic map of degree $d \geq 2$. Assume K is a totally invariant set. Let C denote the critical set of $f$. Assume $\overline{\bigcup_{n=1}^{\infty} f^n(C)} \bigcap K = \emptyset$. Is f hyperbolic on K?
\ep

\bp
Let $H_d$ denote the space of holomorphic maps  $f:\PP^2 \rightarrow \PP^2$ of degree d. This is a finite dimensional space parametrized by the coefficients. Does the set of $ f\in H_d$ with infinitely many attractive basins have measure zero?
\ep

\bp
Let  $f:\PP^2 \rightarrow \PP^2$ be a holomorphic map of degree $d \geq 2$ and let $\lambda_1 \leq \lambda_2$ be the Lyapunov exponents for the ergodic measure $\mu = T \wedge T$. Is $\lambda_1 >0$?
\ep

\bp
Classify critically finite maps on $\PP^3$.
\ep

\newpage

\section{ Complex Analytic Varieties}
\subsection{Embedding Varieties}
\bp Is every compact reduced analytic space biholomorphically equivalent to a
subvariety of a complex manifold? If so, can the manifold be chosen to be compact?\\
-E. L. Stout, stout@math.washington.edu
\ep
\bp  Let C be a smooth cubic in $\PP^2$. Then $\PP^2\setminus C$
is  algebraically parallellizable, but cannot be algebraically immersed in $\CC^2$. Can it be  holomorhically immersed in  $\CC^2$?
[ There are smooth cubics for which it is true . For related problems see the section called ``Stein manifolds''.]\\
-R. Narasimhan, Dept. of Math., Univ. of Chicago, Chicago, IL 60637 
\ep

\subsection{Varieties over the unit disc}
\bp Consider the following statement: `` Let Y be a compact subset of $\CC^2$
which lies over the circle $|\lambda|= 1$, and let U be an open
neighborhood of K, the polynomially convex hull of Y. Then there exists a finitely
sheeted Riemann surface S lying over $|\lambda| < 1$ as a branched covering such
that $S\subset U$.'' Is this statement true for every choice of Y and U?\\
-H. Alexander, u22330@uicvm.bitnet and J. Wermer, Brown University.
\ep
\subsection{Complex Differential Geometry}
\bp Let M be a K$\ddot{a}$hler manifold. Assume that all curvatures are between
$-a^2$ and $-b^2$. Construct nontrivial bounded holomorphic functions on M.\\
-N. Sibony, sibony@anh.matups.fr
\ep

\newpage

\section{ CR manifolds}
\subsection{Embedding}
\bp Are strongly pseudoconvex CR surfaces of real dimension 5 locally embeddable?
\ep
\subsection{Extensions}
\bp Let $M \subset P^n$ be a compact odd dimensional CR manifold. Find conditions
on M which ensure that M bounds a complex variety.(This is well understood in
$\CC^n$.)\\
-N.Sibony, sibony@anh.matups.fr
\ep
\subsection{ Quasiconformal mappings}
\bp {\bf (New) (Liouville type assertion on CR 3-manifolds) On a CR 3-manifold
one can always define a Carnot-Caratheodory metric. A homeomorphism between
two strongly pseudoconvex CR 3-manifolds is called conformal if it maps
infinitesimal spheres with respect to a Carnot-Caratheodory metric to
infinitesimal spheres. Is a conformal homeomorphism CR? \\
-Puqi Tang, tang@math.purdue.edu}
\ep
\bp  {\bf (New) A diffeomorphism between two strongly pseudoconvex CR
$(2n+1)$-manifolds is called quasiconformal if its differential preserves
the underlying contact structures and distorts the CR structures boundedly.
When Hermitian metrics on the contact bundles are fixed, this distortion
can be measured by checking how spheres in the contact space is mapped to
ellipsoids. However, this measurement depends on the choices of the Hermitian
metrics if $n>1$. Fixing a quasiconformal diffeomorphism, can we choose 
Hermitian metrics so that the distortion is minimal?\\
- Puqi Tang, tang@math.purdue.edu}
\ep

\newpage

\section{Scratchpad}

\bp Let $\Omega \subseteq \CC^n$ be a bounded domain with complete Kobayashi
metric.
Is there a constant $C = C(p) > 0$ such that if $f$ is holomorphic and $p > 0$
then for all $z \in \Omega$ and all $r > 0$ small it holds that
$$
|f(z)|^p \leq C \frac{1}{|B(z,r)|} \int_{B(z,r)} |f(\zeta)|^p dV(\zeta) .
$$
Here $f$ is any holomorphic function and $B(z,R)$ is the Kobayashi
metric
ball.  The measure being used is up for grabs.  It can be Euclidean or
one of the canonical measures associated with the Kobayashi metric
construction
(such as that due to Eisenman or Bun Wong).  In case $\Omega$ is finite
type in $\CC^2$ or strongly pseudoconvex in any dimension then the
results is true just because metric balls are comparable to polydiscs
and then classical arguments of Stein (or see my SCV book) will suffice.
In general I think that something like this should be true but I have no
idea how to prove it.\\
-S. Krantz, sk@math.wustl.edu
\ep

\newpage

\section{ Open prize problems}

\bop Prove that if $\Omega$ is strongly pseudoconvex then there is an
absolute constant $K > 0$ such that the integrated Kobayashi distance
between any two points can be realized by a Kobayashi chain of $K$
discs.  [For instance, on a convex domain, $K = 1.$] \\
Prize \$50 -S. Krantz. sk @math.wustl.edu
\eop
\bop Prove that if $\Omega$ is a smooth bounded pseudoconvex domain in $\CC^n$
 which is finite type in the sense of D'Angelo, then it is of finite
 type in the sense of Kohn (ideal type).\\
 Prize \$50 -J. McNeal, mcneal@math.princeton.edu
\eop

\end{document}